\documentclass[10pt,journal]{IEEEtran}

\usepackage{amsmath,cite,xspace}
\usepackage{graphicx}
\usepackage{epsfig}
\usepackage{epstopdf}
\usepackage{bm}
\newcommand{\subparagraph}{}
\usepackage{color}
\usepackage{amssymb}
\usepackage{hyperref}
\newcommand{\utwi}[1]{\mbox{\boldmath $ #1$}}

\newcommand*{\QEDA}{\hfill\ensuremath{\blacksquare}}%


\newenvironment{remark}[1][Remark]{\begin{trivlist}
\item[\hskip \labelsep {\bfseries #1}]}{\end{trivlist}}

\usepackage{mathtools}

\DeclarePairedDelimiter\floor{\lfloor}{\rfloor}

\newcommand{\reals}{\mathbb{R}}

\newcommand{\E}{{\mathbb{E}}}

\newcommand{\diag}{{\textrm{diag}}}
\newcommand{\bc}{{\bf c}}
\newcommand{\ba}{{\bf a}}

\newcommand{\bg}{{\bf g}}
\newcommand{\bh}{{\bf h}}
\newcommand*\bl{\ensuremath{\boldsymbol\ell}}

\newcommand{\bp}{{\bf p}}
\newcommand{\bq}{{\bf q}}

\newcommand{\bx}{{\bf x}}

\newcommand{\bv}{{\bf v}}
\newcommand{\bw}{{\bf w}}

\newcommand{\bz}{{\bf z}}
\newcommand{\by}{{\bf y}}

\newcommand{\bB}{{\bf B}}

\newcommand{\bE}{{\bf E}}

\newcommand{\bR}{{\bf R}}

\newcommand{\bI}{{\bf I}}

\newcommand{\balpha}{{\utwi{\alpha}}}

\newcommand{\bdelta}{{\utwi{\delta}}}

\newcommand{\bphi}{{\utwi{\phi}}}

\newcommand{\bxi}{{\utwi{\xi}}}

\newcommand{\bSigma}{{\utwi{\Sigma}}}

\newcommand{\bmu}{{\utwi{\mu}}}

\newcommand{\cN}{{\cal N}}

\newcommand{\cR}{{\cal R}}

\newcommand{\cB}{{\cal B}}

\newcommand{\cX}{{\cal X}}
\newcommand{\cY}{{\cal Y}}

\newcommand{\sfT}{\textsf{T}}
\renewcommand{\Pr}{\mathbb{P}}
\newcommand{\matE}{\mathbb{E}}

\newcommand{\rev}[1]{{\color{black} #1}}

\graphicspath{{./Figures/}} \DeclareGraphicsExtensions{.eps,.pdf,.jpg,.mps,.png}
\begin{document}

\title{Joint Chance Constraints in AC Optimal Power Flow: Improving Bounds through Learning} 
\author{Kyri Baker, \emph{Member}, IEEE\vspace{-0.8cm} and Andrey Bernstein, \emph{Member}, IEEE}
\maketitle

\begin{abstract}
This paper considers distribution systems with a high penetration of  distributed, renewable generation and addresses the problem of incorporating the associated uncertainty into the optimal operation of these networks. 
Joint chance constraints, which satisfy multiple constraints simultaneously with a prescribed probability, are one way to incorporate uncertainty across sets of constraints, leading to a chance-constrained optimal power flow problem. 
Departing from the computationally-heavy scenario-based approaches or approximations that transform the joint constraint into conservative deterministic constraints, this paper develops a scalable, data-driven approach which learns operational trends in a power network, eliminates zero-probability events (e.g., inactive constraints), and accurately and efficiently approximates bounds on the joint chance constraint iteratively. In particular, the proposed framework improves upon the classic methods based on the union bound (or Boole's inequality) by generating a much less conservative set of single chance constraints that also guarantees the satisfaction of the original joint constraint. The proposed framework is evaluated numerically using the IEEE 37-node test feeder, focusing on the problem of voltage regulation in distribution grids.
\end{abstract}

\section{Introduction} \label{sec:intro}
The AC optimal power flow (OPF) problem is one of the fundamental problems in power system operation and analysis; see, e.g., \cite{OPFoverview} for an overview. \textcolor{black}{Recently, the introduction of dynamic pricing, renewable generation, energy storage, and other distributed energy resources has increased the uncertainty in achieving predictable and reliable grid operations when using deterministic methods. Broadly, there are two classes of approaches which incorporate uncertainties into OPF problems: Robust methods and stochastic methods (and combinations thereof) \cite{Louca16, Muhlpfordt16, Jabr13, Molzahn18}. Amongst stochastic methods, the chance-constrained AC OPF (CC-AC-OPF) is most valuable in situations where the uncertainty lies within constraints, and constraint violations are allowed with a small probability. In power and energy applications, this can include violations in the thermal limit of transmission lines \cite{Vrak13}, relaxation of voltage regulation requirements \cite{DallAnese17CC,Guo18_1}, or loosening of indoor thermal comfort limits \cite{Garifi18}}.

\noindent A prototypical CC-AC-OPF is given by
\begin{subequations} \label{eqn:OPF}
\begin{align} 
 \, \,  &\min_{\substack{\bx \in \cX}}\hspace{.2cm} \matE_{\by} \, f(\bx, \by)  \\
&\mathrm{subject\,to:~} \, \by = \bh(\bx, \bxi) \label{eqn:PF_constr}\\
&\hspace{2cm} \Pr \{\by \in \cY \} \geq 1 -\epsilon\label{eqn:state_constr},
 \end{align}
\end{subequations}
where $\bx$ is a vector that collects all the controllable inputs to the system, typically active and reactive power injections of the controllable distributed energy resources (DERs); $\bxi$ is a random vector representing the uncertainty in the system (e.g., power injections of the uncontrollable assets and solar irradiance); $\by$ is the vector of state variables, such as voltage phasors across the buses of the network;  \eqref{eqn:PF_constr} are the power-flow constraints;  \eqref{eqn:state_constr} is an operational constraint formulated as a chance constraint on the state vector $\by$; and we use the notation $\matE_{\by}$ to denote the expected value with respect to the distribution of $\by$. In particular, \eqref{eqn:state_constr} ensures that the state vector lies in some prescribed operational set $\cY$ with probability at least $1 - \epsilon$ for some (small) quantity $\epsilon > 0$. 

In many applications, the constraint $\by \in \cY$ is composed of several individual constraints $\by \in \cY_i$, $i = 1, \ldots, n$, that have to be satisfied simultaneously;  therefore chance constraint \eqref{eqn:state_constr} is of the form
\begin{equation} \label{eqn:joint_cc}
\Pr \left(\cap_{i=1}^n \{\by \in \cY_i \} \right) \geq 1 -\epsilon.
\end{equation}
Examples of individual constraints that have to be satisfied simultaneously with high probability include joint constraints over different buses in the network,  constraints that link timesteps (e.g., ensuring that power delivered to a sensitive resource is satisfied with high probability across the timesteps after a contingency), or even simply two-sided  constraints (e.g., constraining the upper and lower limits on uncertain line flows or voltage magnitudes).

\textcolor{black}{Many recent works regarding chance constraints in OPF problems have focused on conservative convex upper bounds of single chance constraints \cite{Tyler15,Nemirovski} and distribution-free, data-based single chance constraints \cite{RoaldCC18, Guo18_1, Baker16NAPS, Misra18}. These techniques provide tractable approaches to addressing single chance constraints, and can be used in conjunction with the technique to reduce the joint chance constraint presented in this paper. In addition, robust optimization techniques can be combined with chance-constrained approaches \cite{Venzke18, Lubin16}, which can also be used in conjunction with the formulation in this paper.} \textcolor{black}{Due to the difficulty in handling joint chance constraints, most of these works focus on single chance constraints.} Considering simultaneous probabilistic constraints generally requires either computationally heavy sampling-based approaches which are limited by problem size \cite{Hong_JCC}; or assumptions about the random parameters \cite{JCC06}; or the use of the union bound, or Boole's inequality  \cite{Boole1854}, to \textcolor{black}{separate the joint chance constraint into single constraints and} create conservative upper bounds on the single constraints  \cite{Grosso14, Blackmore09}. In \cite{Hong_JCC}, a Monte Carlo method is proposed to solve a sequence of convex optimization problems, avoiding the use of Boole's inequality, with a guarantee that the algorithm converges to a KKT point. However, it is limited by problem size to small or medium size problems with less than 100 dimensions. Scenario approaches can be used to simplify joint constraints into deterministic single constraints; however, these approaches 
can be overly conservative, and can actually perform worse as the number of samples increases \cite{Campi}. 

Using the union bound (or, Boole's inequality) is the most popular way to relax \eqref{eqn:joint_cc} that boils down to replacing it with $n$ chance constraints
\begin{equation} \label{eqn:single_cc}
\Pr \{\by \in \cY_i \}  \geq 1 -\epsilon_i, \, \, i = 1, \ldots, n.
\end{equation}
It is easy to see that if $\sum_{i = 1}^n \epsilon_i = \epsilon$, \eqref{eqn:single_cc} implies \eqref{eqn:joint_cc}; particularly, if no additional information is used regarding the individual constraints, the typical choice is $\epsilon_i \equiv \frac{\epsilon}{n}$. However, this choice may result in highly conservative solution to \eqref{eqn:OPF}. To illustrate this fact, consider two constraints: $\by \in \cY_1$ and $\by \in \cY_2$. Suppose that the events $A_i := \{\by \notin \cY_i\}$ are highly correlated, in the sense that with very high probability, whenever $A_1$ happens, $A_2$ happens as well (and vice versa). For example, $A_i$ can represent a violation of voltage upper bound at bus $i$ equipped with a photovoltaic (PV) panel, and both buses are geographically close to one another. 
In this case, 
\begin{equation}
\Pr (A_1 \cup A_2 ) = \Pr(A_1) + \Pr(A_2) - \Pr(A_1 \cap A_2) \approx \Pr(A_1) \approx \Pr(A_2)
\end{equation}
because $\Pr(A_1 \cap A_2) \approx \Pr(A_1) \approx \Pr(A_2)$. Therefore, the joint chance constraint \eqref{eqn:joint_cc} would boil down to a single constraint
\begin{align*}
\Pr ( \{\by \in \cY_1\} \cap \{\by \in \cY_2\}) &= 1 - \Pr (A_1 \cup A_2 ) \\
&\approx 1 - \Pr(A_1) \geq 1 - \epsilon,
\end{align*}
or equivalently, $\Pr(A_1) \leq \epsilon$.
However, the union bound approximation  \eqref{eqn:single_cc} will impose a pair of constraints $\Pr(A_i) \leq \frac{\epsilon}{2}$, $i = 1, 2$, therefore unnecessarily restricting the constraint set.

In this paper, we leverage statistical learning tools to address the problem of computationally burdensome joint chance constraints in AC OPF problems, with the following key ingredients:
\begin{itemize}
\item We present a framework for reducing a joint chance constraint into a series of single chance constraints in a method that significantly reduces the conservativeness compared to using Boole's inequality \cite{Baker_boole}. \rev{To this end, we leverage}  support vector classifiers to \emph{classify events $A_i := \{\by \notin \cY_i\}$ as having either zero or non-zero probabilities}. \rev{We term the events (and corresponding constraints) that has non-zero probability as \emph{active}; otherwise, they are \emph{inactive}. That is:
\begin{align*}
    \Pr( \{\by \notin \cY_i\} ) = 0 & \, \Longleftrightarrow \,  \{\by \in \cY_i\} \, \textrm{ is inactive} \\ 
    \Pr( \{\by \notin \cY_i\} ) > 0 & \, \Longleftrightarrow \,  \{\by \in \cY_i\} \, \textrm{ is active} 
\end{align*}
For example, voltage constraints are classified as active or inactive.
}
\item An estimation method is presented which iteratively provides a tighter upper bound on the joint chance constraint and can be terminated before the estimation is finished in computationally restrictive or high dimensional settings where the entire joint constraint cannot be estimated.
\end{itemize}
Unlike classic Monte-Carlo-based approaches, the proposed framework is scalable to high-dimensional constraints. Moreover, the reduction of the joint chance constraint into single chance constraints allows for the use of many of the distributionally-robust single chance constraint reformulations in the literature \cite{Tyler15, Nemirovski}. \textcolor{black}{It is important to note here that we are not developing a new technique for evaluating or reformulating single chance constraints; we are developing a technique for reducing a joint chance constraint into single ones.} Building upon our previous initial work \cite{BakerJointGlobalSIP}, the proposed method can also reduce computation time in non-stochastic settings by removing non-binding constraints from the deterministic optimization problem. 

Simulation results are presented for the IEEE 37-node test system with a high penetration of distributed solar in an active distribution network. While the results presented here are focused on voltage regulation in distribution networks, the method proposed in this paper can be applied in general CC-AC-OPF settings for any type of joint chance constraints.

The remainder of the paper is structured as follows. Section \ref{sec:JCC} discusses the joint chance constraints formulation and outlines our approach. Section \ref{sec:SVC} presents a method to classify inactive constraints and to estimate the remaining joint constraints. Section \ref{sec:systemmodel} outlines the distribution system model and related notation. Section \ref{sec:optimization} discusses the application of the proposed method to voltage regulation problem in active distribution networks. Section \ref{sec:simulations} presents the numerical results. Finally, Section \ref{sec:conclusion} concludes the paper.

\section{Outline of the Approach} \label{sec:JCC}
To explain how we will use statistical learning to reduce the complexity of the joint chance constraint in power network optimization, consider \eqref{eqn:single_cc} and let $A_i := \{\by \notin \cY_i\}$. Then, $\Pr \left(\cap_{i=1}^n \{\by \in \cY_i \} \right) = 1 - \Pr\Big(\bigcup_{i=1}^n A_i\Big)$, and from the probabilistic version of the inclusion-exclusion principle \textcolor{black}{we have the following:}

\begin{align}\label{inclusion}
\begin{split}
\Pr\Big(\bigcup_{i=1}^n A_i\Big) &= \sum_{i=1}^n \Pr(A_i) - \sum_{i<j} \Pr(A_i \cap A_j) + \cdots \\
&\cdots + (-1)^{n-1} \Pr\Big( \bigcap_{i=1}^n A_i \Big)\\
& := \sum_{i=1}^n \Pr(A_i) - P_c.
\end{split}
\end{align}

\begin{figure}[t!]
    \includegraphics[width=0.5\textwidth]{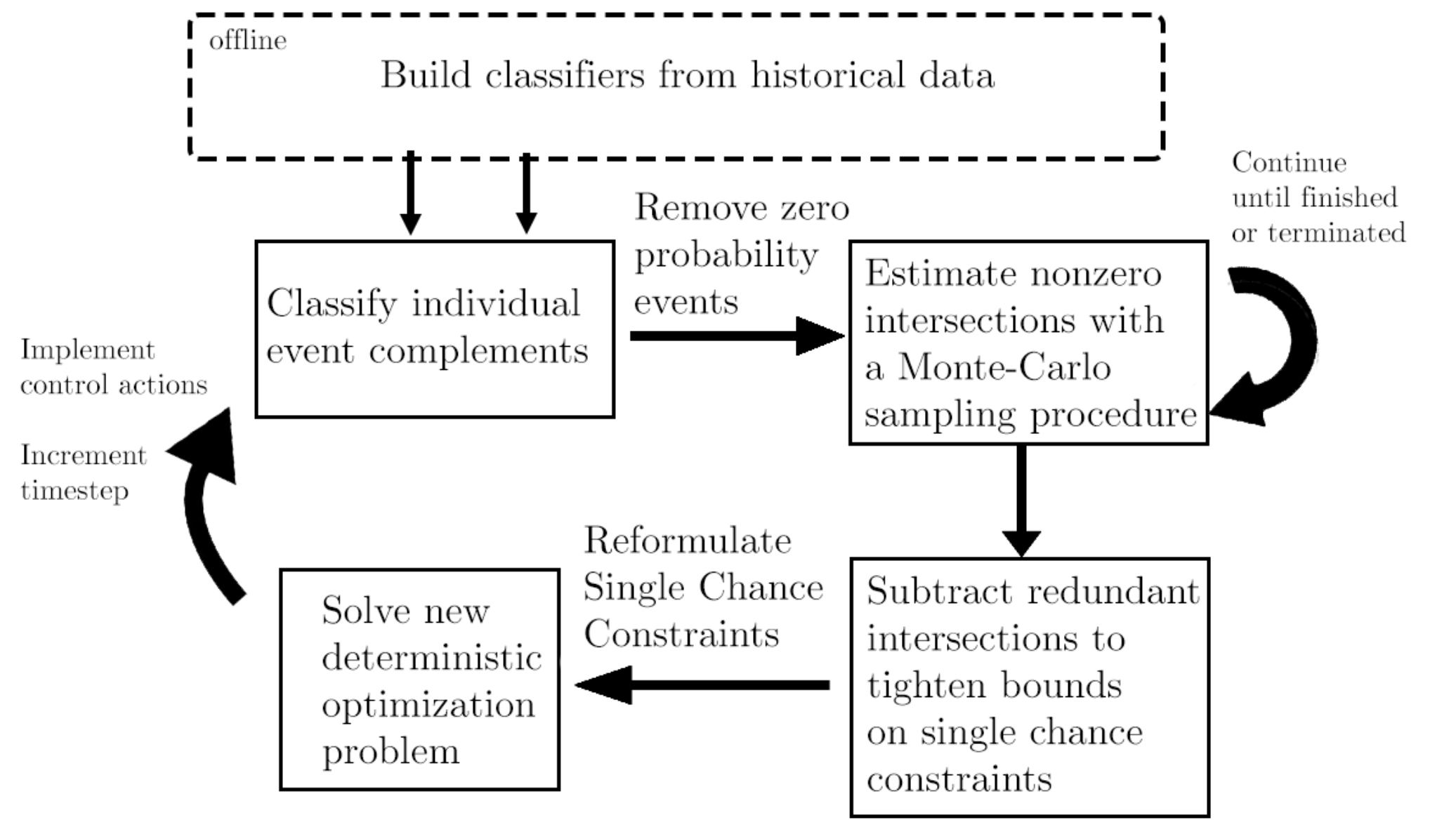}
    \small \caption{An outline of the general procedure for solving the joint chance constrained problem. \textcolor{black}{The details of each individual component in the flowchart can be found in the following sections.}}
    \label{fig:flow}
\end{figure}

\noindent Truncation of the above to Boole's inequality $\Pr(A_1 \cup A_2 \cup ... \cup A_n) \leq \sum_{i=1}^n \Pr(A_i)$ allows for the separation of the joint chance constraint \eqref{eqn:joint_cc} into individual constraints $\Pr(A_i) \leq \epsilon_i$ where $\sum_{i=1}^n \epsilon_i = \epsilon$, \textcolor{black}{and a common choice for $\epsilon_i = \frac{\epsilon}{n}$, where $n$ is the number of individual constraints}. While useful and a very popular technique for solving joint chance constrained programs, Boole's inequality tends to result in very conservative solutions \cite{Baker_boole}. 
To address the deficiencies of using conservative upper bounds, computationally heavy scenario-based approaches, or making assumptions about the single or joint probability distributions, we will present a general, distribution-agnostic technique based on learning marginal probabilities to eliminate zero-probability joint events in \eqref{inclusion} and leverage a Monte Carlo sampling-based approach to estimate the remaining joint probabilities $P_c$. Then, we decompose the joint constraint into single chance constraints that must be satisfied with probability $\geq 1 - \frac{\epsilon}{n} - \frac{P_c}{n}$, \textcolor{black}{where $n$ in this case is the number of nonzero-probability individual events}. Thus, if $P_c > 0$ (i.e., events $A_1, ... , A_n$ are not disjoint), a tighter upper bound for the single chance constraints is provided compared to using Boole's inequality. A general flow chart of the overall procedure is shown in Fig. \ref{fig:flow}; the individual blocks are discussed next.

\section{Constraint Classification and Estimation} \label{sec:SVC}

In an optimization problem, \emph{inactive} constraints are those which, if removed from the problem, would not change the optimal solution. Active constraints, on the other hand, are essential in determining the optimal solution and would change the optimal solution if removed. Machine learning approaches to solve OPF problems have recently been realized as a powerful tool \cite{Misra18, Dobbe18}; here, we leverage machine learning for identifying active constraints in AC OPF problems with joint chance constraints. This section discusses how we can learn which constraints are likely to be inactive in power system optimization given certain system conditions, reducing the computational burden of calculating each term in \eqref{inclusion}.

\subsection{A Simple Example - Two Sided Constraint}
To illustrate the overall idea of the framework, consider the two-sided joint chance constraint which constrains the state of charge $E^{(t+1)}$ of an energy storage system (ESS) to be within desired bounds $\underline{E}$ and $\overline{E}$ with probability at least $1-\epsilon$:

\begin{align} \label{ESS_JCC}
\Pr(\underline{E} \leq E^{(t+1)} \leq \overline{E}) \geq 1 - \epsilon,
\end{align}

\noindent While maintaining ESS state of charge within certain bounds can extend the lifetime of the ESS, under certain situations it may be more beneficial or unavoidable to violate these limits. Intuitively, in certain situations it can be obvious if $\underline{E} \leq E^{(t+1)}$ or $E^{(t+1)} \leq \overline{E}$ is an inactive constraint; for example, if the ESS is currently at its maximum charge ($E^{(t)} = \overline{E}$), and the maximum discharge rate makes it impossible for the ESS to reach $\underline{E}$ in the next time step, we know with certainty that $\underline{E} < E^{(t+1)}$ and thus $\Pr(\underline{E} > E^{(t+1)}) = 0$. So, from the inclusion-exclusion principle, $\Pr(\underline{E} \leq E^{(t+1)} \leq \overline{E}) =1 -  \Pr(\underline{E} > E^{(t+1)}) + \Pr(E^{(t+1)} > \overline{E}) + \Pr(\{\underline{E} > E^{(t+1)}\} \cap \{E^{(t+1)} > \overline{E}\}) = \Pr(E^{(t+1)} \leq \overline{E})$, reducing the joint chance constraint to the single chance constraint $\Pr(E^{(t+1)} \leq \overline{E}) \geq 1 - \epsilon$. However, when dealing with multi-time step problems, assuming one of these events has a zero probability may not be trivial; it may also not be trivial depending on maximum charge/discharge rates, time in between control decisions, the level of uncertainty, or distance between $\underline{E}$ and $\overline{E}$. In addition, while we may have physical intuition as to when a constraint is likely to be relevant or not, there can be many factors influencing the outcome of an optimization problem, and we would like to have an automated way of reducing the complexity of joint chance constraints. Thus, it is desirable to develop a rule that may allow us to exploit these patterns by learning them over time and having the optimization problem automatically decompose the joint chance constraints into single chance constraints depending on the outcome of these rules.

In general, recall that if $\Pr(A_i) = 0$, $\Pr(A_i \cap A_j) = 0$ for all $A_j$; if even a single constraint is classified as inactive, a significant number of terms in the joint chance constraint expansion are eliminated from the calculations and do not have to be estimated further. As a larger example, consider a four-event union $\Pr(A_1 \cup A_2 \cup A_3 \cup A_4)$ and its expansion via \eqref{inclusion}:

\small
\begin{align}
\begin{split}
&\Pr(A_1) + \Pr(A_2) + \Pr(A_3) + \Pr(A_4)\\ 
&- \Pr(A_1 \cap A_2) - \Pr(A_1 \cap A_3) - \Pr(A_1 \cap A_4) \\
&- \Pr(A_2 \cap A_3) - \Pr(A_2 \cap A_4) - \Pr(A_3 \cap A_4) \\
&+ \Pr(A_1 \cap A_2 \cap A_3) + \Pr(A_1 \cap A_2 \cap A_4) \\
&+ \Pr(A_1 \cap A_3 \cap A_4) + \Pr(A_2 \cap A_3 \cap A_4)\\
&- \Pr(A_1 \cap A_2 \cap A_3 \cap A_4)
\end{split}
\end{align}
\normalsize

\noindent If constraint $A_1$ is classified as inactive, the above reduces to

\small
\begin{align}
\begin{split}
&\Pr(A_2) + \Pr(A_3) + \Pr(A_4)\\ 
&- \Pr(A_2 \cap A_3) - \Pr(A_2 \cap A_4) - \Pr(A_3 \cap A_4) \\
&+ \Pr(A_2 \cap A_3 \cap A_4)
\end{split}
\end{align}
\normalsize

\noindent dramatically reducing the number of intersections we must estimate. For sizable joint chance constraints, identifying zero probability events can potentially make an otherwise intractable problem possible to solve via sampling approaches.

\begin{figure*}[!htb]
\begin{center}
    \includegraphics[width=0.8\textwidth]{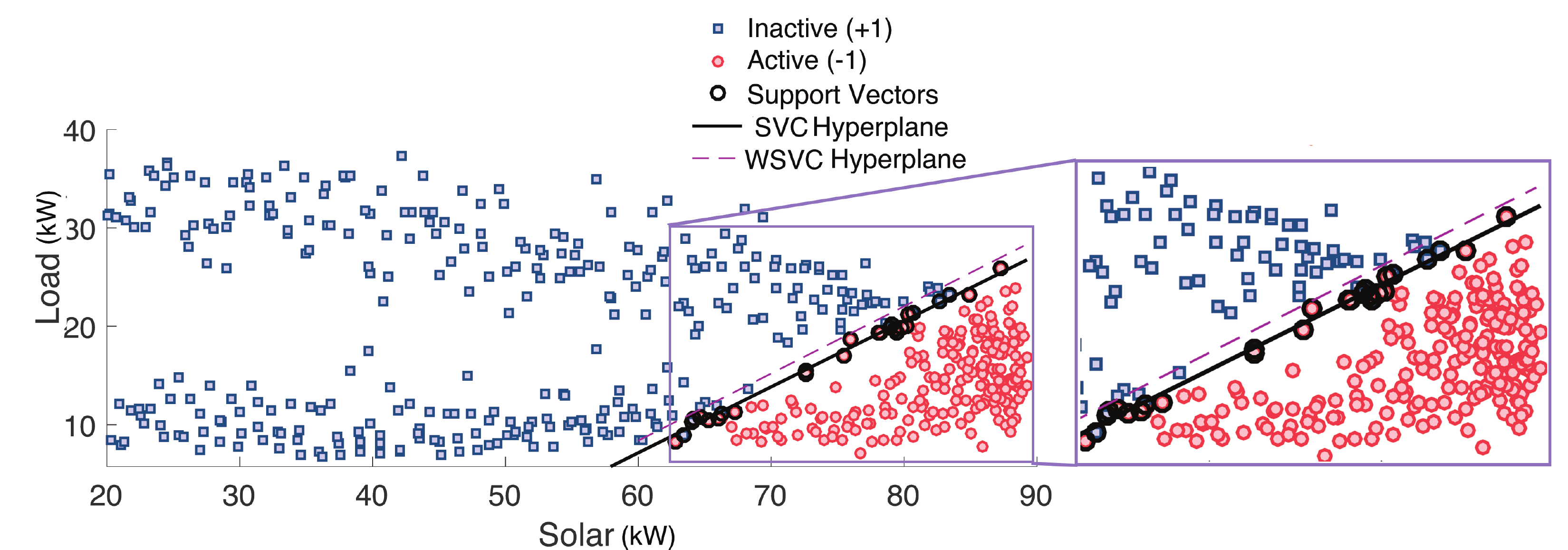}
    \small \caption{Visualizing the support vector classifier for classifying overvoltage events ($V > \overline{V}$) at a node in a distribution network. Intuitively, we know that overvoltages occur when solar exceeds load by a particular amount, and the classifiers provide us with a selection rule for including or excluding voltage constraints. Left: Resulting support vector classifiers from traditional SVC and a weighted version (wSVC) which heavily penalizes any misclassifications of active constraints as inactive. Right: A magnified version. The support vectors are training data which lie closest to the separating hyperplane and, if removed, would change the solution of (P0).}
    \label{fig:SVC}
  \end{center}
\end{figure*} 

\subsection{Support Vector Classification (SVC)}

\textcolor{black}{Next, we will develop classifiers for classifying constraints as active or inactive. This procedure is performed before the final OPF problem is solved to reduce the dimensionality of the joint chance constraint.} We will use a popular machine learning technique for classification called Support Vector Classification (SVC). \textcolor{black}{For each training sample $i = 1,...,m$,} and two classes, namely active ($\ell_i = -1$) and inactive ($\ell_i = +1$), we wish to create a \textcolor{black}{decision} rule which uses selected inputs to determine whether or not we include that constraint in the optimization problem. Here, we seek to form an affine classifier of the form $\bw^T\bphi + b$ with weights $\bw \in \mathbb{R}^2$ and bias $b \in \mathbb{R}$ that classifies constraints as active  ($\bw^T\bphi + b \geq 0$) or inactive  ($\bw^T\bphi + b < 0$) based on input features $\bphi \in \mathbb{R}^2$ (e.g., load and available solar at a node in the distribution network as considered in the example in Section \ref{sec:simulations} below). In our formulation, called ``weighted SVC (wSVC)", we heavily penalize misclassifications of active constraints as inactive, while maximizing the separation between classes \cite{StatLearn}. \textcolor{black}{Unlike typical applications for support vector machines/classifiers where misclassifications are equally weighted, highly weighting misclassifications of active constraints as inactive pursues the preservation of an upper bound on the original joint constraint. If some inactive constraints were classified as active, the bound may get looser but stay valid - the other way around, and the bound may not be preserved.} In the training stage, we build the classifier by using $m$ samples of labeled training data $\bl$, by solving the following optimization problem:

\begin{subequations} 
\label{SVC_opt}
\begin{align} 
 \mathrm{(P0)}  &\min_{\substack{\bw, \bz, b}}  \hspace{.2cm}  \frac{1}{2}\bw^Tb + \bc^T\bz \label{SVC_obj}\\
 &\mathrm{s.t.} ~~~\ell_i(\bw^T\bphi_i + b) \geq 1 - z_i\label{SVC_c}\\
&~~~~~~~z_i \geq 0 \label{SVC_z}
\end{align}
\end{subequations}

\noindent where $\bc \in \mathbb{R}^{m}$ is a penalty parameter and $c_i = 1$ if $\ell_i = +1$ and $c_i = a$ if $\ell_i = -1$, $a \gg 1$. \textcolor{black}{Using this objective weighting for $\bz$ in \eqref{SVC_obj}, we deviate from the traditional SVC formulation with the wSVC problem, which heavily penalizes misclassifications of active constraints as inactive.} Each slack variable $z_i$ is nonzero if $\bphi_i$ is classified incorrectly, and zero otherwise. \textcolor{black}{After solving (P0), a classifier $\bw^T\bphi + b$ for each individual node is formed, which takes the form of an affine function of the two features: the total load at a node and available solar generation at a node.}

\subsection{Classifying Power System Constraints} \label{sec:class}
\rev{We next describe how we classify constraints for a particular power-system problem. For illustration purposes, we focus on 
the example of distribution grid voltage regulation under high PV production; the corresponding OPF problem (P1) is formally defined in Section \ref{sec:optimization} below. At each node in the considered distribution network, we construct an SVC/wSVC with the training data set consisting of:
\begin{itemize}
    \item Features: net solar production and load at that node.
    \item Labels: $\ell = 1$ if the voltage magnitude at that node is below the uppder limit; $\ell = -1$ otherwise.
\end{itemize}
The features in the training data set are constructed from historical data on solar and load, whereas the labels are computed by solving the OPF (P1) \emph{without voltage constraints}, and verifying whether the corresponding voltage limits are violated or not.
As overvoltage conditions are primarily caused when solar generation exceeds consumption (load),  it is reasonable to exploit this relationship to use the current and forecasted levels of solar and load to determine which voltage constraints will be relevant or binding. The classifiers are then constructed by solving (P0) and obtaining corresponding affine decision rules.

In Fig \ref{fig:SVC}, we illustrate the differences between the classifier chosen with the traditional SVC versus the conservative weighted SVC (wSVC). 
}
Traditional SVC and wSVC may not perfectly separate the two classes, but the conservative wSVC \textcolor{black}{aims to ensure} that all training points labeled as active are correctly classified. 
\rev{
\begin{remark}[Remark 1]
In the voltage regulation application and particular formulation used in this paper, the classes can be divided by a separating hyperplane. In applications with a nonlinear relationship between classes, additional methods to preprocess/transform the data can be employed. For example, lifting to a higher dimensional space whereby the data is linearly separable can be achieved, e.g., using kernel-based methods \cite{StatLearn}.
\end{remark}

\begin{remark}[Remark 2]
While in this paper we focus on using these classifiers for probabilistic constraints, the approach would also provide computational benefits for constraint removal from deterministic programs as well. In particular with the voltage regulation case - as seen in Section \ref{sec:simulations}, including upper bound constraints on the voltage is unnecessary throughout most of the daily operating period.
\end{remark}
}

\subsection{Iteratively Estimating Event Intersections} \label{sec:iterative}
The number of intersections given in the joint chance constraint expansion that must be estimated is given by

\begin{align} \label{numterms}
\sum_{k=2}^{|M|} {{|M|}\choose k} 
\end{align}

\noindent where $M$ is the set of indices of active constraints, $|M|$ is the cardinality of $M$ (i.e., the number of active constraints), and $n$ is the total number of constraints. For example, if the joint chance constraint originally contained 8 constraints, 247 intersections must be estimated to recover the original constraint. If half of these were classified as inactive, only 11 intersections must be estimated, which is much more reasonable for solving optimization problems on fast timescales. The expansion of the joint chance constraint can now be written as

\begin{align}\label{inclusion2}
\begin{split}
\Pr\Big(\bigcup_{i=1}^n A_i\Big) &= \sum_{m\in M} \Pr(A_m) - \sum_{m \neq j, m, j \in M} \Pr(A_m \cap A_j) + \cdots \\
&\cdots + (-1)^{|M|-1} \Pr\Big( \bigcap_{m \in M} A_m \Big).
\end{split}
\end{align}

\noindent 

\begin{figure}[t!]
    \includegraphics[width=0.5\textwidth]{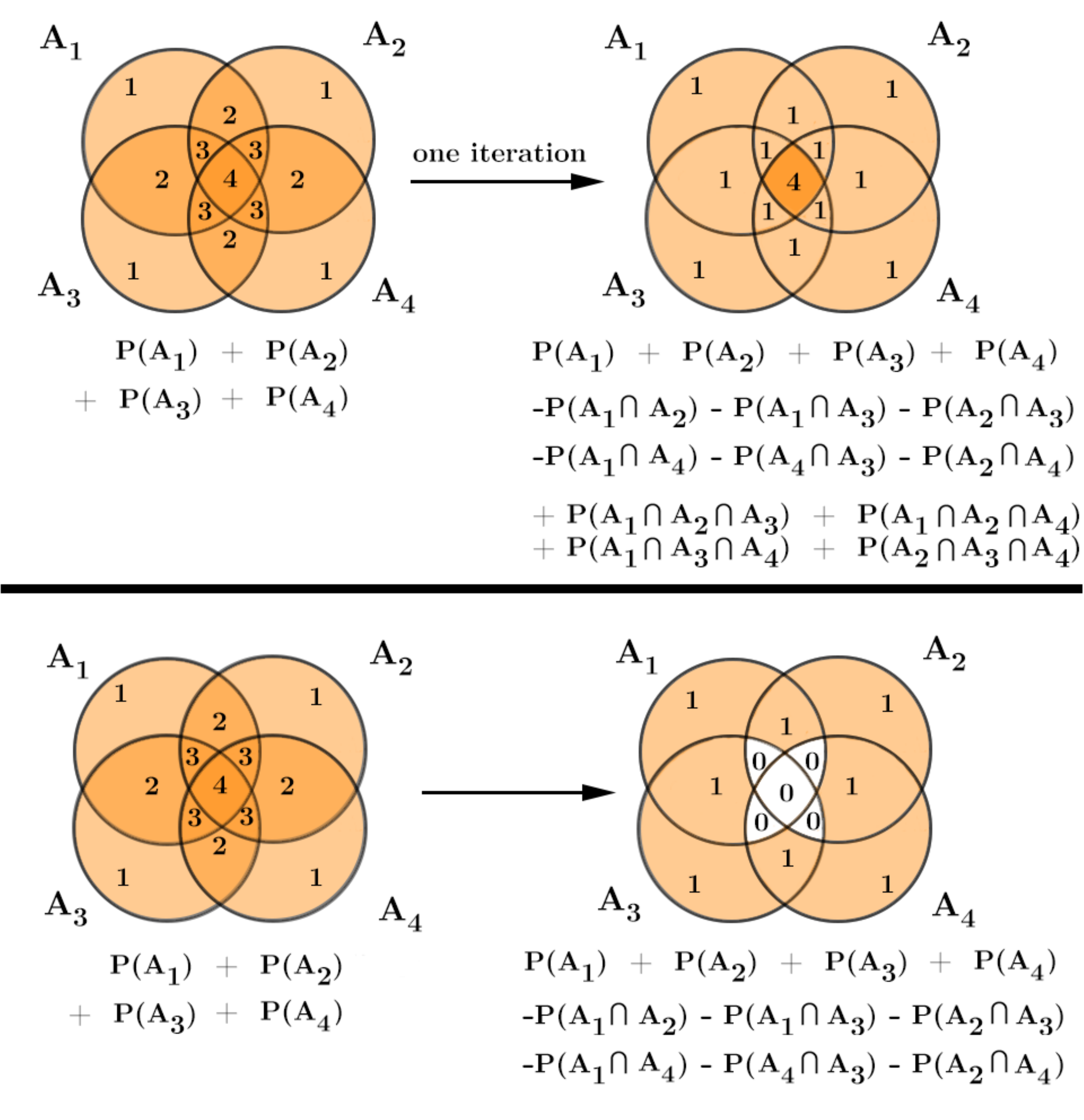}
    \small \caption{A four-event example where the numbers overlaid on the events and intersections represent the number of times that intersection is accounted for. The union bound (left) tends to overestimate the intersection of events, resulting in overly conservative control actions. Top: An improved upper bound on $\Pr(A_1 \cup A_2 \cup A_3 \cup A_4)$ is sought by removing redundant intersections (right) as time, data availability, and problem size allow, maintaining an upper bound on the original constraint by performing pairwise intersection estimations (here, only one iteration of the intersection estimation algorithm is performed). Bottom: If the intersections are simply accounted for in the order that \eqref{inclusion2} provides, an upper bound may not be preserved.}
    \label{fig:Boole_tot}
\end{figure}

Our goal is to iteratively estimate event intersections in a way that maintains an upper bound on the original joint chance constraint, allowing for the termination of the algorithm before the entire joint chance constraint is estimated. 
In fact, the order in which these intersections are computed is very important; if certain intersections are included in the expansion but not others, an upper bound of the original union of events may not be preserved. 
We therefore only estimate intersection probabilities for \emph{pairs} of terms in \eqref{inclusion2}:

\begin{align}\label{newbound}
\begin{split}
\cB_K :=& \sum_{m \in M} \Pr(A_m) \\
& -\sum_{k=1}^K \Big[\sum_{\substack{I \subset {\{1,...,|M|\}} \\ |I| = 2k}} \Pr\Big(A_I\Big) - \sum_{\substack{I \subset {\{1,...,|M|\}} \\ |I| = 2k+1}} \Pr\Big(A_I\Big)\Big]
\end{split}
\end{align}

\noindent for $K = 1 ... \floor{\frac{|M|+1}{2}}$, where $I \subset \{1,...,|M|\}, |I| = k$ denotes all subsets $I$ of indices $1, ..., |M|$ which contain exactly $k$ elements, and $A_I := \bigcap_{i \in I} A_i$. A four-event example shown in Fig. \ref{fig:Boole_tot} to illustrate this: in the top subfigure, an improved upper bound on $\Pr(A_1 \cup A_2 \cup A_3 \cup A_4)$ is sought by removing redundant intersections (right) as time, data availability, and problem size allow, maintaining an upper bound on the original constraint by performing pairwise intersection estimations (here, only one iteration of the intersection estimation algorithm is performed). In the bottom subfigure, the event intersection probabilities are iteratively removed in the order of \eqref{inclusion2}, no longer maintaining an upper bound on the union of events. \textcolor{black}{In this example, the calculation in \eqref{newbound} only performs one iteration at $k = 1$; the combinations of pairwise intersections is calculated for the sum where $|I| = 2k = 2$, and the combinations of three-way intersections is calculated for the sum where $|I| = 2k+1 = 3$. These two terms are then added to the marginal probabilities $\Pr(A_m)$.}

This provides a benefit over the convenient but extremely conservative Boole's inequality as well as a more reliable and robust alternative to scenario-based approaches, which may require more time than available in between control actions in large networks. In the worst case (no computation time is allowed to estimate intersections), the algorithm is equivalent to using Boole's inequality to create tractable single chance constraints.  

\begin{remark}[Observation 1] \label{proofbound}
We have that
\[
\Pr\Big(\bigcup_{i=1}^n A_i\Big) \leq \cB_K \leq \sum_{m \in M} \Pr(A_m)
\]
for all \vspace{1.5mm} $K \in \{1, \ldots, \floor{\frac{|M|+1}{2}}\}$.
\end{remark}

\noindent \emph{Proof.} 
The proof follows by the inclusion-exclusion principle, the monotonicity of $\cB_K$ in $K$, and the fact that for $K = \floor{\frac{|M|+1}{2}}$, $\Pr\Big(\bigcup_{i=1}^n A_i\Big) = \cB_K$.

\QEDA

Note that Observation 1 allows us to terminate the estimation process of the joint probabilities before $K = \floor{\frac{|M|+1}{2}}$. \textcolor{black}{This process is particularly useful when the number of non-zero terms in the joint chance constraint is large, as we can ensure that an upper bound is preserved which is still tighter than that provided by Boole's inequality.}

In this paper, we estimate the remaining joint probabilities \rev{(i.e., the second term in \eqref{newbound}) using a sampling approach, and represent these probabilities with their relative frequencies. For that purpose, a deterministic optimization problem (e.g., (P1) in Section \ref{sec:optimization}) is solved in a similar way as described in Section \ref{sec:class} using historical inputs (solar, load), and the relative frequency of the event intersections are computed from that data. For example, if an event $A_1 \cap A_2$ occurred 3,000 times out of 10,000, we would assign $\Pr(A_1 \cap A_2) = \frac{3,000}{10,000} = 0.3$. This process is discussed in more detail in Section \ref{sec:estimate}.}


\section{Distribution Network and System Models}
\label{sec:systemmodel}

Consider a distribution feeder comprising $N$ $PQ$ nodes and a single slack node. Let $V_n \in \mathbb{R}$ 
denote the  line-to-ground voltage magnitude
at node $n$, 
and define the $N$-dimensional vector  $\bv := [V_1, \ldots, V_N]^\sfT \in \mathbb{R}^{N}$.
Constants $P_{\ell,n}$ and $Q_{\ell,n}$ denote the real and reactive demands at node $n$, and we can define the vectors  $\bp_\ell := [P_{\ell,1}, \ldots, P_{\ell,N}]^\sfT$ and $\bq_\ell := [Q_{\ell,1}, \ldots, Q_{\ell,N}]^\sfT$; if no load is present at node $n$, then $P_{\ell,n} = Q_{\ell,n} = 0$.

Here, we use a linearization of the AC power-flow equations \cite{sairaj2015linear,linModels} which linearly relates the voltage magnitudes $\bv$ to the injected real and reactive powers $\bp \in \reals^N$ and $\bq \in \reals^N$ in the form
\begin{align} 
\bv & \approx \bR \bp + \bB \bq + \ba, \label{eq:approximate} 
\end{align}

\noindent where $\bR$, $\bB$, and $\ba$ are parameters that are dependent on the system model \cite{sairaj2015linear,linModels} \textcolor{black}{and are usually dependent on system line parameters, topology, and substation voltage, but can also be computed from data-driven techniques such as regression-based methods \cite{BakerNetwork17}}. While the proposed methodology does not require problem convexity, we leverage a linearization in order to provide a clear exposition of the joint chance constraint reformulation in Section \ref{sec:optimization}.

\subsection*{Photovoltaic (PV) Systems}
Random quantity $P_{av,n}$ denotes the maximum renewable-based generation at node $n$  -- hereafter referred to as the available solar power. Particularly, $P_{av,n}$ coincide with the maximum power point at the AC side of the inverter. When RESs operate at unity power factor and inject the available solar power $P_{av,n}$, issues related to power quality and reliability in distribution systems may be encountered. For example, voltages exceeding prescribed limits at a particular node may be experienced when RES generation exceeds the load of that consumer~\cite{PVhandbook}.  Efforts to ensure reliable operation of existing distribution systems with increased behind-the-meter renewable generation are focus on the possibility of inverters providing reactive power compensation and/or curtailing real power. To account for the ability of the RES inverters to adjust the output of real power, let $\alpha_n \in [0,1]$ denote the fraction of available solar power curtailed by RES-inverter $n$. If no PV system/inverter is at a particular node $i$, $P_{av,i} = \alpha_i = 0$. For convenience, define the vectors $\balpha :=[\alpha_1,\ldots, \alpha_N  ]^\sfT$ and $\bp_{av} := [P_{av,1}, \ldots, P_{av,N}]^\sfT$. 

The available active power generation from solar is modeled as $\bp_{av} = \overline{\bp}_{av} + \bdelta_{av}$, where $\overline{\bp}_{av} \in \mathbb{R}^N$ is a vector of the forecasted values and $\bdelta_{av}  \in \cR_{av}  \subseteq \mathbb{R}^{N}$ is a random vector whose distribution captures spatial dependencies among forecasting errors. We assume that the distribution system operator has a certain amount of information about the probability distributions of the forecasting errors $\bdelta_{av}$. This information can come in the form of either knowledge of the  probability density functions, or a model of $\bdelta_{av}$ from which one can draw samples. In this paper, we make the assumption that these errors are normally distributed \textcolor{black}{with zero mean; i.e., $\bdelta_{av}$ \texttt{{\raise.17ex\hbox{$\scriptstyle\mathtt{\sim}$}}} $\mathcal{N}(0, \sigma)$, and thus the remaining single chance constraints can be exactly reformulated as analytical expressions and included in the optimization problem directly \cite{Roald13}.} However, distributionally robust formulations of single chance constraints \cite{Tyler15, Nemirovski} can easily be incorporated into the framework here, \textcolor{black}{and thus the Gaussian assumption for the random quantities is not necessary for our framework}.

\section{Joint Chance Constrained Formulation}
\label{sec:optimization}
\subsection{Optimization problem reformulation}
The joint chance constraint optimization for voltage regulation in distribution systems is shown below:

\begin{subequations} 
\label{PJCC}
\begin{align} 
 \mathrm{(P1)}&  \min_{\substack{\bv, {\tiny \balpha} }}  \hspace{.2cm}  \E(f(\bv, \balpha, \bp_{\ell}, \bq_{\ell})) \label{eq:Pmg_cost}\\
& \mathrm{subject\,to} \nonumber \\
& \bv  = \bR ((\bI - \diag\{\balpha\}) \bp_{\mathrm{av}} -  \bp_{\ell}) \nonumber \\
& \hspace{1cm}- \bB \bq_{\ell} + \ba\label{mg-balance-t} \\
& \Pr\{V_{1} \leq V_{\mathrm{max}}, ... , V_{N} \leq V_{\mathrm{max}} \} \geq 1- \epsilon \label{mg-prob-jcc}\\
& 0 \leq \alpha_i \leq 1, \, i = 1, \ldots, N; \hspace{3.05cm}  \label{mg-alpha} 
\end{align}
\end{subequations}
\noindent Constraint~\eqref{mg-balance-t} represents a surrogate for the power balance equation; constraint \eqref{mg-prob-jcc} is the joint chance constraint that requires voltage magnitudes to be within $V_{max}$ with at least $1-\epsilon$ probability; and constraint \eqref{mg-alpha} limits the curtailment percentage from $0-100\%$. The cost function $f(\bv, \balpha, \bp_{\ell}, \bq_{\ell})$ is convex and can consider a sum of penalties on curtailment, penalties on power drawn from the substation, penalties on voltage violations, etc.

By reformulating the joint constraint (\ref{mg-prob-jcc}) as a series of single chance constraints, we can write

\begin{align}
\Pr(V_{i} \leq V_{\mathrm{max}} ) \geq 1- \epsilon_i \label{singlecc}
\end{align}

\noindent for all $i \in 1, \ldots, N$, and each $\epsilon_i$ is chosen such that $\sum_{i=1}^N \epsilon_i = \epsilon$. In the case of Boole's inequality, we choose $\epsilon_i = \frac{\epsilon}{N}$; in the case of what we call the \emph{improved} Boole's inequality, using the method proposed in this paper, we choose $\epsilon_i = \frac{\epsilon}{N} + \frac{P_c}{N}$, where $P_c$ is our estimation of the non-zero probabilities in \eqref{inclusion}. If $P_c > 0$ (all events are not mutually exclusive), it is clear that the improved Boole's inequality provides a less conservative upper bound on the chance constraints.

\subsection{Estimating Event Intersections}\label{sec:estimate}
We use a conservative relative frequency sampling approach to estimate the event intersections that have not been classified as zero. If event $A_i$ represents an overvoltage at node $i$, and $I$ is a subset of nodes, we can estimate the probability of the intersection of overvoltages at the nodes included in $I$ as 

\begin{align} \label{MC}
\Pr\left(\bigcap_{i \in I} A_i\right) \approx \frac{\sum_{s=1}^{N_s} \textbf{1} \left \{\bv_I(\bdelta_s) > V_{max}\right\} }{N_s}~~,
\end{align}

\noindent for $N_s$ random draws of the uncertainty distribution, where draw $s$ is denoted $\bdelta_s$. Vector $\bv_I$ contains the voltage magnitudes at each of the nodes in $I$, and $\textbf{1} \left \{\bv_I(\bdelta_s) > V_{max}\right\}$ is one if \emph{all} of the elements in $\bv_I(\bdelta_s)$ are greater than $V_{max}$ and zero otherwise. To represent the most conservative case, for each sample $\bdelta_s$, the voltage vector $\bv(\bdelta_s)  = \bR (\bI(\overline{\bp}_{\mathrm{av}} + \bdelta_s -  \bp_{\ell}))- \bB \bq_{\ell} + \ba$; i.e., the curtailment variables $\balpha$ are chosen to be zero to represent no curtailment and thus the most conservative case for the control policy. 
\rev{Note that, using the objective function defined in Section \ref{sec:simulations} below, this is equivalent to solving a deterministic version of problem (P1) without voltage constraints.}
The impact of different sample sizes $N_s$ and computational burden of the estimation process is discussed in further detail in the next section.

\subsection{Analytical Reformulation of Single Chance Constraints}

The single chance constraints can be reformulated as exact, tractable constraints \cite{BoVa04}, assuming $\epsilon \leq 0.5$. Assuming the joint distribution of the random variables is a multivariate Gaussian with  mean $\bmu$ and covariance matrix $\bSigma$, define $\mu_i$ as the $i$-th value in $\bmu$ and $\sigma_i$ as the $(i,i)$-th entry in $\bSigma$. Then, define the following function at each node $i$:

\begin{align*}
h(p_{av,i}) = &\sum_j(R_{ij}[(1-\alpha_j)p_{av,j} - p_{\ell, j}]) \\
&- \sum_j(B_{ij}q_{l,j}) + a_i - V_{max}
\end{align*}

\noindent where $R_{ij}$ is the $(i,j)$-th entry of $\bR$, $B_{ij}$ is the $(i,j)$-th entry of $\bB$, and $a_i$ is the $i$th element of $\ba$. Then $h(p_{av,i})$ is also normally distributed with the following mean $\mu'_i$ and variance $\sigma'_i$:

\begin{align*}
\mu'_i &= \sum_j(R_{ij}[(1-\alpha_j)\mu_j - p_{\ell, j}]) - \sum_j(B_{ij}q_{l,j}) + a_i - V_{max}\\
\sigma'_i &=  \sum_jR_{ij}(1-\alpha_j)\sigma_j
\end{align*}

\noindent Thus, the constraints \eqref{singlecc} can be reformulated using the Gaussian cumulative distribution function (CDF) $\Phi$:

\begin{align*}
\Pr\{h(p_{av,i}) \leq 0\} = \Phi\Big(\frac{0-\mu'_i}{\sigma'_i}\Big) \geq 1-\epsilon_i
\end{align*}

\noindent With the final analytical constraint written using the quantile function (the inverse of the Gaussian CDF):

\begin{align}
R_i[(1-\alpha_i)\mu_i - p_{\ell,i}] - B_iq_{\ell,i} + a_i - V_{max} \nonumber \\
\leq -R_i\alpha_i\sigma_i\Phi^{-1}(1-\epsilon_i) \label{final_CC}
\end{align}

\noindent Which can be explicitly included into problem (P1) for each $i$ in place of the joint constraint \eqref{mg-prob-jcc}.

\begin{remark}[Remark 3]
 In these results, the individual solar forecasting errors are modeled as Gaussian. Because of this, the single chance constraints can be exactly analytically reformulated. Without loss of generality, other distributionally robust methods for single chance constraints can also be used here \cite{Baker16NAPS, Nemirovski}, but as the contribution of this paper is in the decomposition of the joint chance constraint, not in addressing the tractability of single chance constraints, we have kept the marginal distributions Gaussian for simplicity of exposition. The method proposed in this paper is not distribution-specific.
\end{remark}

\begin{remark}[Remark 4] Note that the original use of Boole's inequality ensures the satisfaction of the original constraint by choosing $\epsilon_i$ such that $\sum_{i=1}^n \epsilon_i = \epsilon$; however, without optimizing this parameter, suboptimal performance of this reformulation is possible \cite{JCC_Opt}. We leave the optimal choice of $\epsilon_i$ as a direction for future work.
\end{remark}

\section{Numerical Results} \label{sec:simulations}

The IEEE-37 node test feeder \cite{IEEE37} was used for the following simulations. Five-minute load and solar irradiance data from weekdays in August 2012 was obtained from \cite{Bank13} for the simulations, and in order to emulate a situation with high-PV penetration and risks of overvoltage, 8 $200$-kW rated PV systems were placed at nodes 29-36. The considered cost function seeks to minimize renewable curtailment; specifically, 

\begin{align} 
f(\bv, \balpha, \bp_{\ell}, \bq_{\ell}) =  \sum_{i \in \cN} d_i  \alpha_i^2 ,  \label{eq:cost_sim}
\end{align} 

\noindent where the cost of curtailing power at each node is set to be $d_i = \$0.10$. The number of samples used to calculate each intersection was $N_s = 10,000$. The considered joint chance constraint considers maintaining voltages at nodes $29 - 36$. Each $\mu_i$, $i=1...N$ was chosen to be the power generated from the forecasted PV at that node, based on the shape of the aggregate solar irradiance from \cite{Bank13} and shifted using samples from a uniform distribution from +/- 1 kW across each node. The covariance matrix $\Sigma$ was formed by setting each entry $(i,j)$ to $\Sigma_{ij} = \bE[(P_{av,i} - \mu_i)(P_{av,j} - \mu_j)^T]$. Three cases were considered in the following numerical results. First, a deterministic case was considered, which uses a certainty equivalence formulation and uses the mean of each of the uncertain parameters in place of each random variable in the optimization problem. Second, Boole's inequality was used to separate the joint chance constraint into a series of conservative single constraints, each with $\sum_{i=1}^N \epsilon_i = \epsilon$ \textcolor{black}{and $\epsilon_i = \frac{\epsilon}{n}$} . Third, an Improved Boole's inequality is considered, where the proposed methodology is implemented to approximate each of the intersections in \eqref{inclusion2}, and $\epsilon_i = \frac{\epsilon}{N} + \frac{P_c}{N}$. 

\subsection{Training, Testing, and Choosing the Number of Samples}
Each of the classifiers (one per constraint; 8 classifiers total) were trained using 1152 samples (4 training days), and tested using 864 samples (3 testing days), using $c_i = 1$ when $\ell_i = +1$ and $c_i = 10$ when $\ell_i = -1$. The overall classification error was $0.19\%$ for false classification of binding events and $4.73\%$ for false classification of non-binding events. A larger classification error for the $\ell_i = +1$ is to be expected from our formulation in (P0); it is more detrimental to the performance of the algorithm if we exclude constraints by mistake rather than include non-binding constraints unnecessarily. \textcolor{black}{If desired, although it has the potential of increasing the conservativeness of the solution, the classifier bias $b$ can also be increased to ensure that the classification error of active constraints as inactive occurs even less frequently. Despite the misclassification of $0.19\%$ of data points in this direction, as the results in the next section demonstrate, the original joint chance constraint is still satisfied within the prescribed probability $1-\epsilon$.}

\begin{figure}[t!]
    \includegraphics[width=0.5\textwidth]{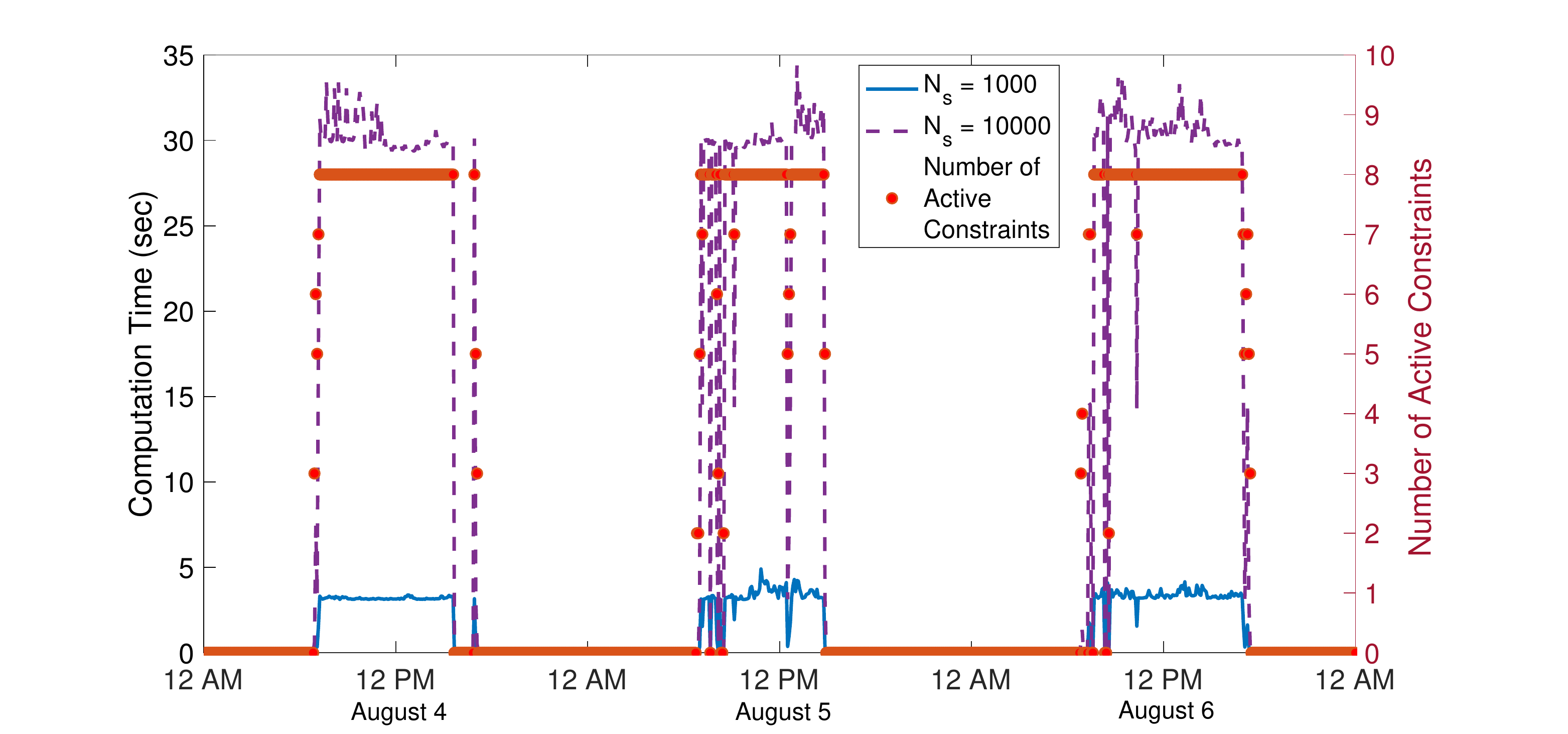}
    \small \caption{The number of active constraints (out of 8 total) is indicated with red dots for each time instance. The total computation time (s) required for calculating the corresponding intersections is shown in black (for $N_s = 1,000$ samples) and purple (for $N_s = 10,000$ samples).}
    \label{fig:comp_time}
\end{figure}

\begin{figure}[t!]
    \includegraphics[width=0.49\textwidth]{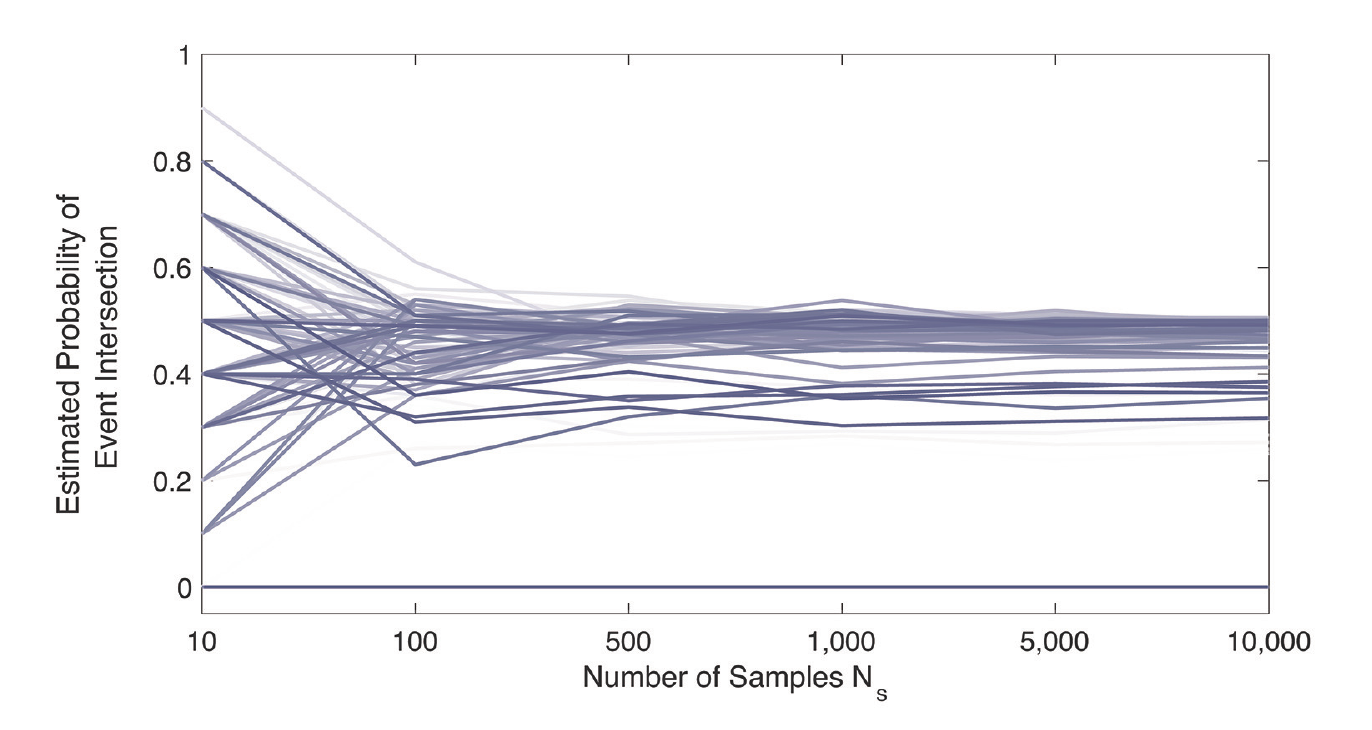}
    \small \caption{A sensitivity analysis showing the estimated probabilities for each of the event intersections in a single test day as a function of the number of samples $N_s$ used to estimate that intersection. The change between estimated intersections with $N_s = 5,000$ and $N_s = 10,000$ is sufficiently small; any additional increase in the number of samples will not provide much additional information.}
    \label{fig:sensitivity}
\end{figure}

In Fig. \ref{fig:comp_time}, the number of active constraints (out of 8 total) is indicated with red dots for each time instance \textcolor{black}{and referenced against the right-hand y-axis}. The total computation time (s) required for calculating the corresponding intersections is shown in black (for $N_s = 1,000$ samples) and purple (for $N_s = 10,000$ samples). Computational time can be reduced by potentially sacrificing accuracy of estimating event intersections; in Fig. \ref{fig:sensitivity}, the value of increasing the number of samples wanes around $N_s = 1,000$. Thus, in the following simulations, the conservative choice of $N_s = 10,000$ was made to estimate each event intersection; however, in a general setting this number is dependent on the underlying distribution.

\subsection{Voltage Regulation Results}

\begin{figure}[t!]
\hspace*{-2mm}    \includegraphics[width=0.49\textwidth]{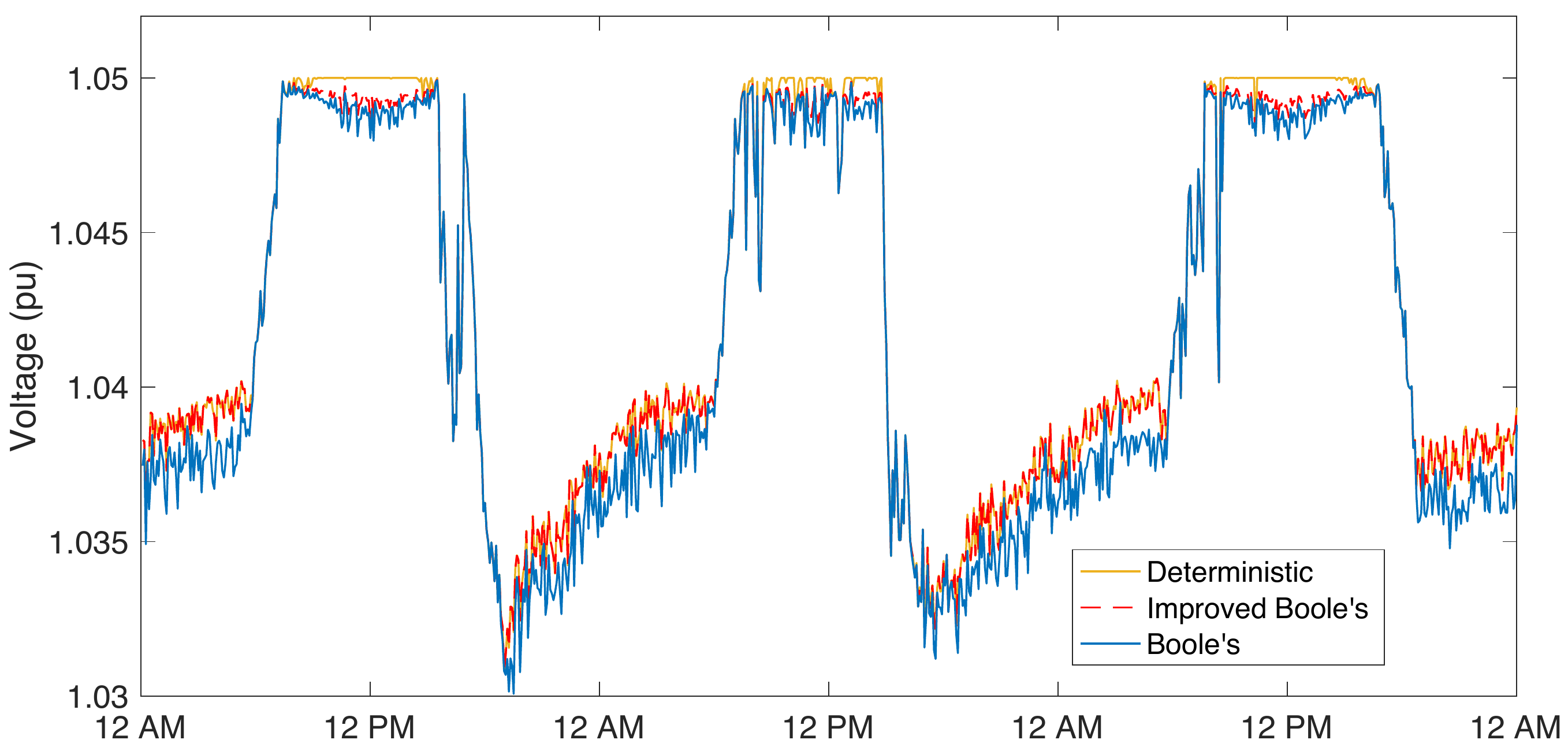}
    \small \caption{The predicted maximum voltages using the control policies determined for voltage regulation in each of the three cases. In the non-deterministic cases, the control policies are more conservative in order to account for the uncertainty in the solar irradiance.}
    \label{fig:control_policy}
\end{figure}

In the following results, the maximum joint constraint violation probability was set to $\epsilon = 0.02$ for the Boole's and Improved Boole's cases, and all of the terms in \eqref{inclusion2} were estimated. In Fig. \ref{fig:control_policy}, the maximum voltage magnitudes from the resulting control policies are shown for each of the three cases. The deterministic case \textcolor{black}{does not} take forecast uncertainty into account, and as a result, the voltages are pushed to the maximum voltage of $1.05$ pu. The Boole's case curtails enough solar generation to ensure that overvoltages will not occur with a high probability; the Improved Boole's case reduces this probability and results in less curtailment.

\begin{figure}[t!]
    \includegraphics[width=0.49\textwidth]{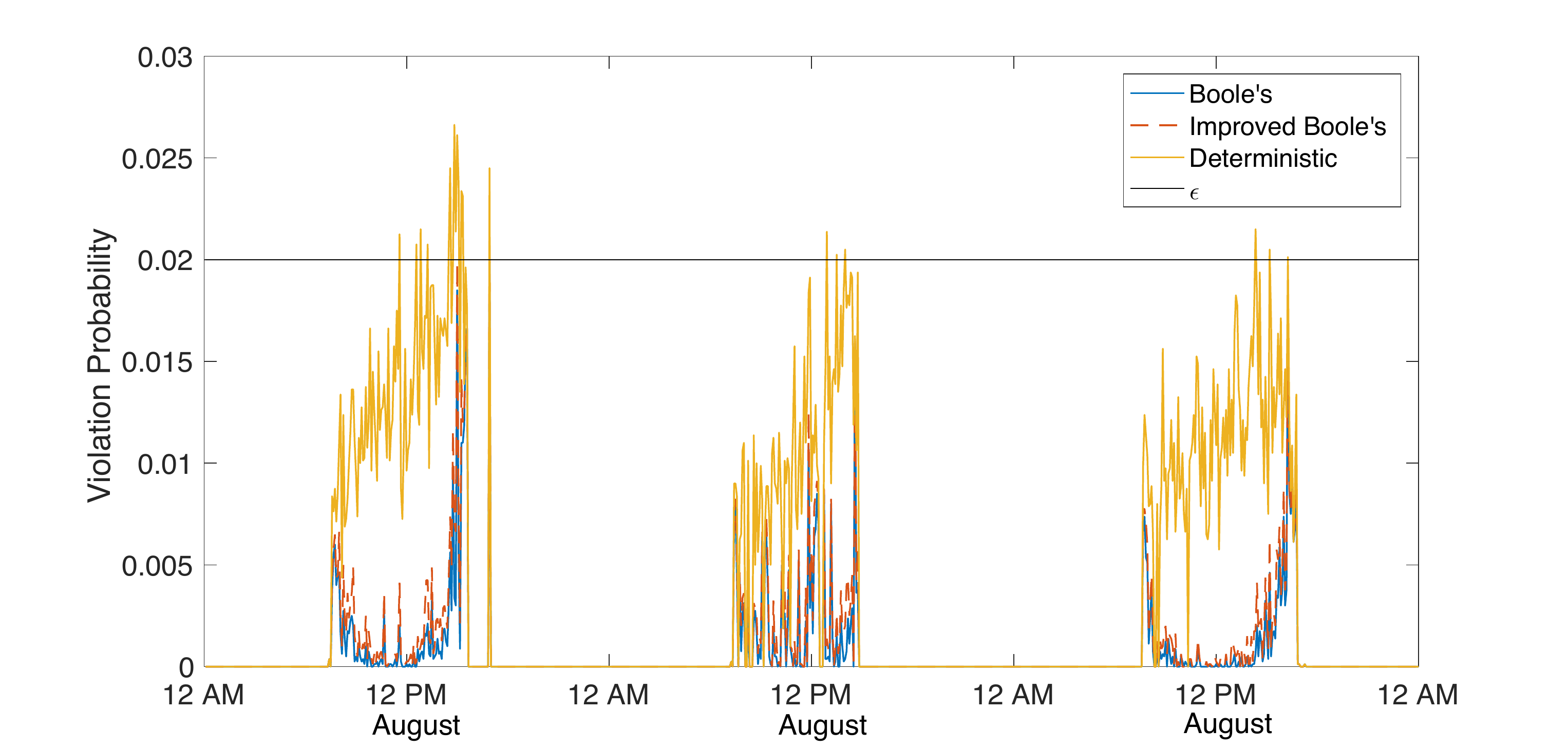}
    \small \caption{The violation probability of the joint chance constraint calculated through a Monte Carlo validation procedure with $N_m = 10,000$. During times of high solar irradiance, the deterministic control policy does not guarantee satisfaction of the joint chance constraint, while the stochastic solutions ensure that the constraint is satisfied with probability $\geq 98\%$.}
    \label{fig:epsilon}
\end{figure}

A Monte Carlo validation procedure was implemented to demonstrate the behavior of the control policies for $N_m = 10,000$ random draws of the uncertainty distributions at each timestep. In Fig. \ref{fig:epsilon}, these resulting probabilities are shown. The deterministic case, which only considers the mean of the random variables, violates the desired chance constraint bound of $0.02$ when compared with the chance constrained methods, because that method offers no guarantee that the voltages will be within limits. Boole's method is generally more conservative than the Improved Boole's method and results in lower violation probabilities, with both methods resulting in satisfaction of the original joint chance constraint.

Table \ref{tab:obj} shows the total objective function value and voltage violation probability across the three-day test period for the deterministic, Boole's, and Improved Boole's cases. As expected, the deterministic case results in the lowest cost (but highest probability of voltage violations); the Boole's case is overly conservative, resulting in a higher level of curtailment and thus cost, but lowest probability of voltage violations. The Improved Boole's case strikes a balance between the two, resulting in a slightly higher violation probability than the original Boole's case but with a lowered objective value.

\begin{table}[t!]
\centering
\caption{Total objective function value and maximum observed voltage violation probability.}
\label{tab:obj}
\begin{tabular}{|l|l|l|l|}
\hline
                                                                          & Deterministic & Boole's & Improved Boole's \\ \hline
\begin{tabular}[c]{@{}l@{}}Total Objective \\ Function Value\end{tabular} & 148.40        & 162.63  & 156.93           \\ \hline
\begin{tabular}[c]{@{}l@{}} Maximum Violation\\ Probability\end{tabular}   & 2.66\%        & 1.85\%  & 1.98\%           \\ \hline
\end{tabular}
\end{table}

\subsection{\textcolor{black}{Computational Time and Multi-Phase Systems}} \label{sec:threephase}
\begin{figure}[t!]
    \includegraphics[width=0.52\textwidth]{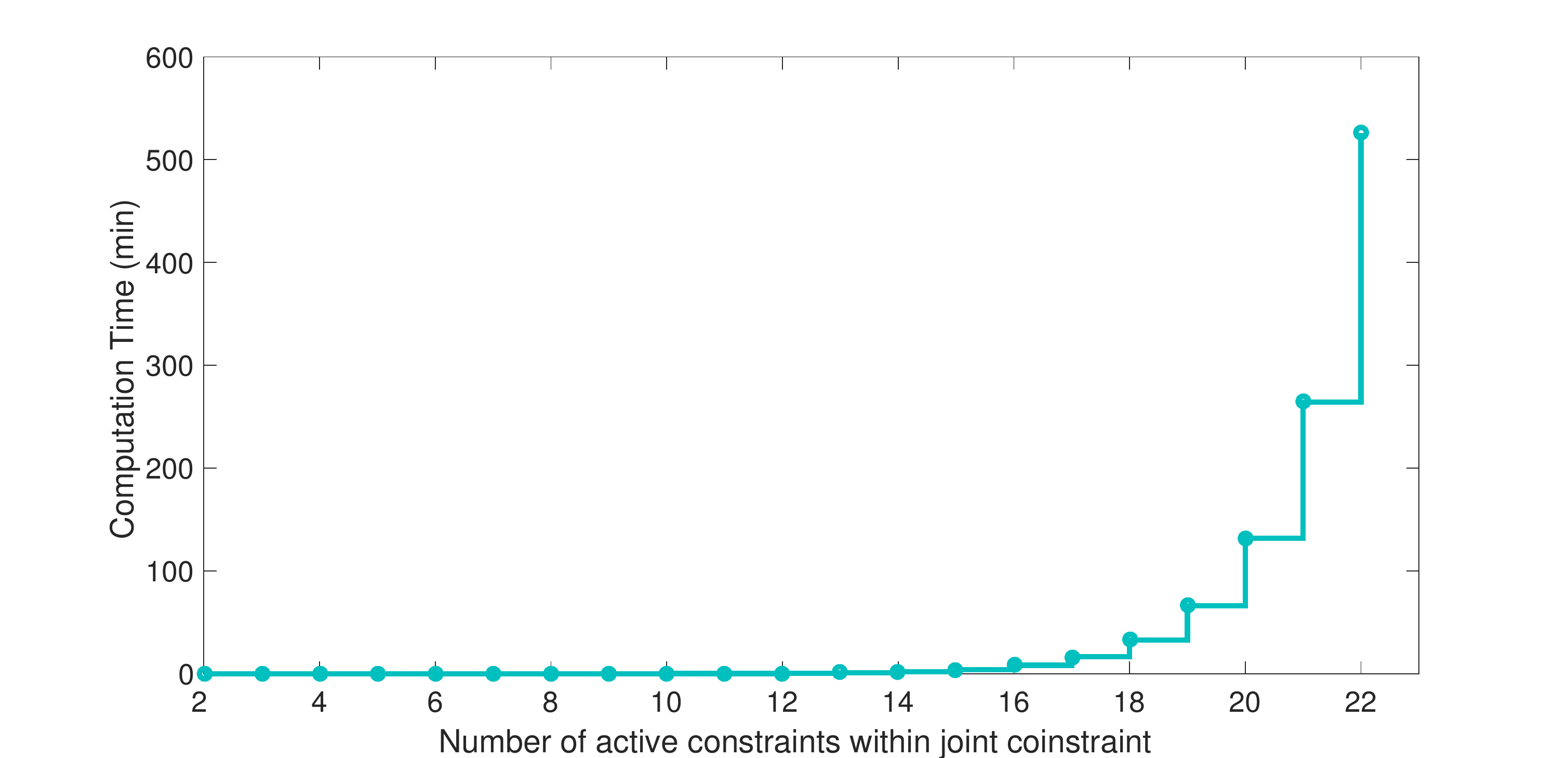}
    \small \caption{\textcolor{black}{The time required to estimate the nonzero terms in the joint chance constraint expansion quickly becomes computationally intractable as $K$ grows. This shows the benefit of the iterative technique proposed in Section \ref{sec:iterative} which allows the algorithm to be terminated prematurely.}}
    \label{fig:comp_time_22}
\end{figure}

\textcolor{black}{To demonstrate the computational burden of the proposed framework as the joint chance constraint increases in number of terms, we have performed additional simulations on an unbalanced, three-phase version of the aforementioned 37-node feeder and using the multi-phase linearization procedure from \cite{linModels}. Three-phase, wye-connected PV systems were connected to the same nodes as in the previous test case, thus increasing the number of single constraints within the joint constraint by threefold (24 total terms) by constraining the voltage magnitude at each of the three phases. As Fig. \ref{fig:comp_time_22} and equation \eqref{numterms} illustrates, estimating the number of intersections for $|M| = 24$ would require more computational time than typically given between OPF control actions. Thus, we can show the benefit of the iterative method proposed in Section \ref{sec:iterative}; for very large joint chance constraints, the algorithm can be terminated prematurely while still providing a tighter or equal upper bound to that of Boole's inequality. This allows the system operator to tune the time required for the constraint estimation procedure to a desired level.}  


\begin{figure}[t!]
    \includegraphics[width=0.52\textwidth]{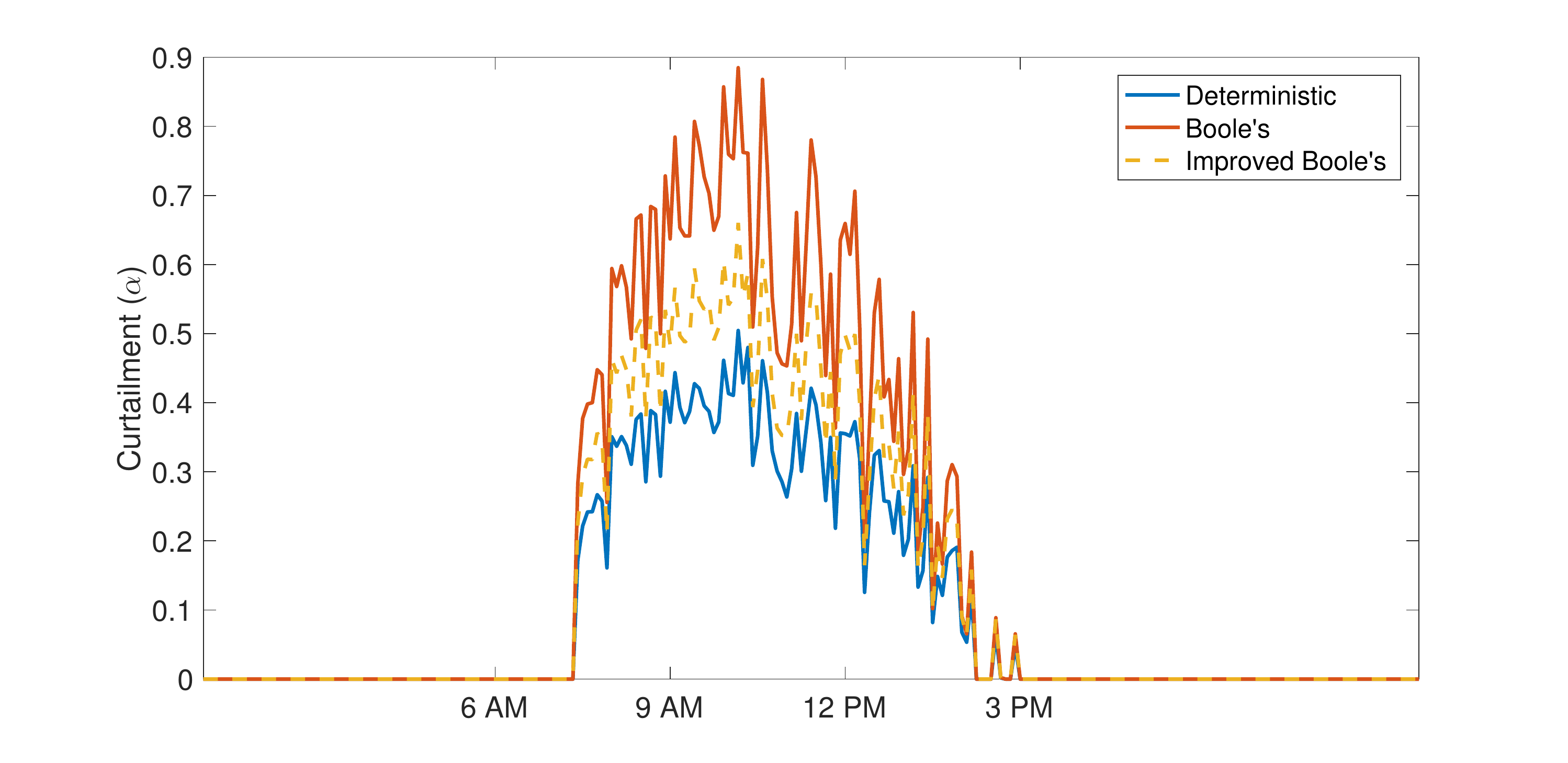}
    \small \caption{\textcolor{black}{The PV curtailment on phase 1 of node 28. As expected, the Improved Boole's method produces a slightly less conservative solution (i.e., less curtailment) than the solution obtained with Boole's inequality.}}
    \label{fig:node28_ph1}
\end{figure} 

\textcolor{black}{For example, here we can terminate the procedure at 12 active constraints, which can be computed in approximately one minute. In this case, over one test day and using $\epsilon = 0.05$ and Monte Carlo simulations with $1,000$ samples per time step, the maximum observed violation probability with Boole's inequality was $3.04\%$ and with the Improved Boole's method was $3.55\%$. While they both provided conservative solutions that were under the prescribed limit of $5\%$, the Improved Boole's method was less conservative. Fig. \ref{fig:node28_ph1} demonstrates this relationship at the inverter at node 28 by showing lower curtailment levels for the solution that used the Improved Boole's method. This tradeoff can also be compared with the computational burden tradeoff - it is likely that if more than 12 active constraints were considered in the constraint estimation procedure, the method would have produced a greater violation probability.}

\section{Conclusion} \label{sec:conclusion}
In this paper, we demonstrated how identifying zero-probability events with support vector classifiers can increase the computational efficiency of computing joint chance constraints via sampling methods. In addition, we provided an iterative approach appropriate for fast timescale optimization, ensuring that if the entire constraint cannot be computed, that the resulting approximation of the constraint always provides an upper bound of the original constraint at every iteration which is tighter than that provided by Boole's inequality. Simulation results were shown which addressed voltage regulation in distribution networks with high PV penetration, and the proposed method was demonstrated to result in a lower cost than Boole's inequality and lower constraint violation probability than a deterministic certainty equivalence formulation.

Future work will address the question of how to optimally allocate the estimated intersection probabilities $P_c$ to the individual chance constraints (rather than allocating them equally across single constraints as in this paper), determining how many samples are adequate for estimating event intersections, and identifying which statistical learning techniques are best suited for identifying active constraints in power systems optimization problems. In addition, an important question for future work is how to incorporate sampling error and uncertainty when identifying and estimating the underlying probability distributions.

\bibliographystyle{IEEEtran}
\bibliography{references,biblioNew}

\begin{IEEEbiography}[{\includegraphics[width=1in]{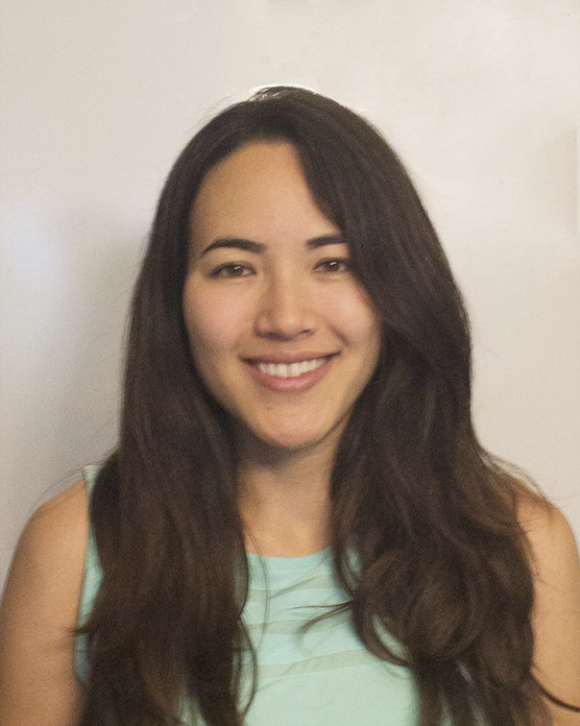}}]{Kyri Baker} (S'08, M'15)
received her B.S., M.S, and Ph.D. in Electrical and Computer Engineering at Carnegie Mellon University in
2009, 2010, and 2014, respectively. Since Fall 2017, she has been an Assistant Professor at the University of Colorado, Boulder, in the Department of Civil, Environmental, and Architectural Engineering, with a courtesy appointment in the Department of Electrical, Computer, and Energy Engineering. Previously, she was a Research Engineer at the National Renewable Energy Laboratory in Golden, CO. Her research interests include power system optimization and planning, building-to-grid integration, smart grid technologies, and renewable energy.
\end{IEEEbiography}

\begin{IEEEbiography}[{\includegraphics[width=1in]{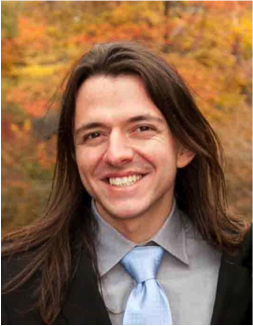}}]{Andrey Bernstein} (M'15)
received the B.Sc. and M.Sc. degrees in Electrical Engineering from the Technion - Israel Institute of Technology in 2002 and 2007 respectively, both summa cum laude. He received the Ph.D. degree in Electrical Engineering from the Technion in 2013. Between 2010 and 2011, he was a visiting researcher at Columbia University. During 2011-2012, he was a visiting Assistant Professor at the Stony Brook University. From 2013 to 2016, he was a postdoctoral researcher at the Laboratory for Communications and Applications of Ecole Polytechnique Federale de Lausanne (EPFL), Switzerland. Since October 2016 he has been a Senior Scientist at NREL. In February 2019, he has become a Group Manager of the Energy Systems Control and Optimization group at NREL.
His research interests are in the decision and control problems in complex environments and related optimization and machine learning methods, with particular application to intelligent power and energy systems. Current research is focused on real-time optimization of power distribution systems with high penetration of and machine learning methods for grid data analytics.
\end{IEEEbiography}


\end{document}